\documentclass[letterpaper, 10 pt, conference]{ieeeconf}  

\IEEEoverridecommandlockouts                              
\usepackage{mathrsfs}
\usepackage{graphics} 
\usepackage{graphicx}
\usepackage{dsfont}
\usepackage{bbm}
\usepackage{amsmath} 
\usepackage{amssymb}  

\usepackage{xcolor}

\let\proof\relax  \usepackage{amsthm}

\newtheorem{theorem}{Theorem}
\newtheorem{proposition}{Proposition}
\newtheorem{lemma}{Lemma}

\theoremstyle{definition}
\newtheorem{definition}{Definition}
\newtheorem{remark}{Remark}

\title{ \LARGE\bf On the robustness of stabilizing feedbacks for  quantum spin-$\bf\frac{1}{2}$  systems
}

\author{Weichao Liang$^{\dag}$, Nina H. Amini$^{*}$ and Paolo Mason$^{*}$
\thanks{$^{\dag}${W. Liang is with Laboratoire Analyse G\'eom\'etrie Mod\'elisation, CY Cergy Paris Universit\'e. 2, av. Adolphe Chauvin, 95302 Cergy-Pontoise, cedex, France.
{\tt\small weichao.liang@u-cergy.fr} }} 
\thanks{$^{*}$ N. H. Amini and P. Mason are with Laboratoire des signaux et syst\`emes (L2S)
CNRS, CentraleSup\'elec, Universit\'e Paris-Sud, Universit\'e Paris-Saclay.
3, rue Joliot-Curie, 91192 Gif-sur-Yvette, cedex, France.
        {\tt\small firstname.lastname@l2s.centralesupelec.fr}}
}

\begin{document}
\maketitle
\thispagestyle{empty}
\pagestyle{empty}
\begin{abstract}
In this paper, we consider stochastic master equations describing the evolution
of quantum spin-$\frac{1}{2}$ systems interacting with electromagnetic fields undergoing continuous-time measurements. We suppose that the initial states and the exact values of the physical parameters are unknown. We prove that the feedback stabilization strategy considered in~\cite{liang2018exponential} is robust to these imperfections. This is shown by studying the asymptotic behavior of the coupled stochastic master equations describing the evolutions of the actual state and the estimated one under appropriate assumptions on the feedback controller. We provide sufficient conditions on the feedback controller and a valid domain of estimated parameters which ensure exponential stabilization of the coupled system. Furthermore, our results allow us to answer positively to~\cite[Conjecture~4.4]{liang2019exponential_TwoUnknown} in the case of spin-$\frac 1 2$ systems with unknown initial states, even in presence of imprecisely known physical parameters. 
\end{abstract}
\section{Introduction}
Classical stochastic filtering theory~\cite{kallianpur2013stochastic,xiong2008introduction} provides 
tools to optimize the estimation of dynamics described by stochastic differential equations, in presence of noisy observations. A primitive theory of quantum filtering theory appeared  in the work of Davies in the 1960s~\cite{davies1969quantum,davies1976quantum}. Belavkin in the 1980s established original results in quantum filtering and feedback control of quantum systems, as a natural extension of classical filtering and control
~\cite{belavkin1983theory,belavkin1989nondemolition,belavkin1992quantum,belavkin1995quantum}. The development of quantum probability theory
and quantum stochastic calculus~\cite{hudson1984quantum,hudson2003introduction,meyer2006quantum} provided essential mathematical tools to describe open quantum systems and quantum filtering. In the physics community, quantum filtering theory is also known as quantum trajectory theory, after it has been established in a more heuristic manner, by Carmichael in the 1990s~\cite{carmichael1993open}. 
A modern introduction to quantum filtering theory may be found in~\cite{bouten2009discrete}. 

Continuous-time quantum filters describe the time
evolution of the states of open quantum systems interacting with electromagnetic fields undergoing continuous-time measurements. Quantum filters are 
solutions of matrix-valued
stochastic differential equations called stochastic master equations.

Quantum filtering theory plays a major role in the development of quantum feedback control. A measurement-based feedback is designed based on the information obtained from quantum filters. The systematic design of stabilizing feedback  control laws for quantum systems is a crucial step towards engineering of quantum devices. In particular, feedback stabilization of pure states has received particular interest~\cite{sayrin2011real,mirrahimi2007stabilizing,van2005feedback}. 
 In real experiments different sorts of imperfections may be present, as for instance inefficient detectors, unknown initial states, imprecise knowledge of the detector efficiency and other physical parameters, etc. Hence, from a practical point of view, choosing feedback controls which are robust to such imperfections is an important, and challenging, problem.

Concerning quantum angular momentum systems with known initial states and parameters,
in~\cite{van2005feedback}, based on numerical approaches, the authors designed for the first time a quantum feedback controller that globally stabilizes a quantum spin-$\frac{1}{2}$ system towards an eigenstate of $\sigma_z$ in presence of imperfect measurements. More recently, in~\cite{mirrahimi2007stabilizing}, by analyzing the stochastic flow and by using stochastic Lyapunov techniques, the authors constructed a switching feedback controller which globally stabilizes the $N$-level quantum angular momentum system, in presence of imperfect measurements, to the target eigenstate. In~\cite{liang2018exponential,liang2019exponential}, by using stochastic and geometric control tools, we provided sufficient conditions on the feedback control law ensuring almost sure exponential convergence to a predetermined eigenstate of the measurement operator for spin-$\frac 12$ and spin-J systems respectively (see~\cite{cardona2020exponential,cardona2018exponential} for exponential stabilization results via a different approach).

In~\cite{liang2019exponential_TwoUnknown}, we considered controlled quantum spin-$\frac12$ systems in the case of unawareness of initial states and in presence of measurement imperfections. We proved that the fidelity between the quantum filter and the associated estimated filter converges to one under appropriate assumptions on the feedback controller. For spin-J systems, we considered feedback controls of a particular form, and we conjectured that such control laws are capable of exponentially stabilize the system towards an eigenstate of the measurement operator.

In this paper, we study the feedback exponential stabilizability of spin-$\frac 12$ systems in presence of measurement imperfections and unawareness of the initial states and of the physical parameters (namely, the detection efficiency, the difference between the energies of the excited state and the ground state, and the strength of the interaction between the system and the probe). We 
find general conditions on the feedback control guaranteeing
robust exponential stabilization with respect to 
such imperfections. 
The dynamics is defined by a coupled system of equations describing the evolutions of the quantum filter and the associated estimated filter, with the feedback controller being a function of the estimated quantum filter.  In order to show our main result, Theorem~\ref{Theorem:Feedback Exp stabilization}, we analyze the asymptotic behavior of the coupled system under appropriate assumptions on the feedback controller. In particular, we provide sufficient conditions on the feedback controller and a valid domain for the estimated parameters which ensure exponential feedback stabilization of the coupled quantum spin-$\frac12$ system. Moreover, we give explicit forms of feedback controllers which guarantee such feedback exponential stabilization. We precise that the stabilizing feedback controllers that we proposed in~\cite{liang2018exponential} satisfy the assumptions of Theorem~\ref{Theorem:Feedback Exp stabilization}
and also that the results of this paper prove~\cite[Conjecture~4.4]{liang2019exponential_TwoUnknown} for spin-$\frac 12$ systems and for a more complicated case, 
since in~\cite{liang2019exponential_TwoUnknown} we assumed unknown initial conditions but precise knowledge of the physical parameters.
Numerical simulations are provided in order to illustrate our results and to support the efficiency of the proposed candidate feedback.
\paragraph{Notations}
The imaginary unit is denoted by $i$. We indicate by $\mathds{1}$ the identity matrix. 
We denote the conjugate transpose of a matrix $A$ by $A^*.$ The function $\mathrm{Tr}(A)$ corresponds to the trace of a square matrix $A.$ The commutator of two square matrices $A$ and $B$ is denoted by $[A,B]:=AB-BA.$

We denote by $\mathrm{int}({S})$ the interior of a subset $S$ of a topological space  and by $\partial{S}$ its boundary. 
\section{System description}
We consider quantum spin-$\frac{1}{2}$ systems. In the following we describe the evolutions of the actual quantum state and its associated estimated state assuming that the initial state and the physical and experimental parameters are not known.  The corresponding coupled system is given by the following stochastic master equations, in It\^o form
\begin{align}
d\rho_t&=L_{\omega, M}^u(\rho_t)dt+G_{\eta,M}(\rho_t)\big(dY_t-\sqrt{\eta M}\mathrm{Tr}(\sigma_z \rho_t)dt\big),\nonumber\\
d\hat{\rho}_t&=L_{\hat{\omega}, \hat M}^u(\hat{\rho}_t)dt+G_{\hat{\eta},\hat{M}}(\hat{\rho}_t)\big(dY_t-\sqrt{\hat{\eta} \hat{M}}\mathrm{Tr}(\sigma_z \hat{\rho}_t)dt\big),\nonumber
\end{align}
where
\begin{itemize}
\item the actual quantum state of the spin-$\frac12$ system is denoted as $\rho$, and belongs to the space $
\mathcal{S}_2:=\{\rho\in\mathbb{C}^{2\times 2}|\,\rho=\rho^*,\mathrm{Tr}(\rho)=1,\rho \geq 0\}$. The associated estimated state is denoted as $\hat{\rho}\in\mathcal{S}_2$,
\item the matrices $\sigma_x,$ $\sigma_y$ and $\sigma_z$ correspond to the Pauli matrices. 
\item $L_{\omega, M}^u(\rho):=-i/2[\omega\sigma_z+u\sigma_y,\rho]+M/4(\sigma_z\rho\sigma_z-\rho)$ and $G_{\eta,M}(\rho):=\sqrt{\eta M}/2\big(\sigma_z\rho+\rho\sigma_z-2\mathrm{Tr}(\sigma_z\rho)\rho\big)$. 
\item $Y_t$ denotes the observation process of the actual quantum spin-$\frac12$ system, which is a continuous semi-martingale whose quadratic variation is given by $\langle Y,Y\rangle_t=t$. Its dynamics satisfies $dY_t=dW_t+\sqrt{\eta M}\mathrm{Tr}(\sigma_z\rho_t)dt$, where $W_t$ is a one-dimensional standard Wiener process,
\item $u:=u(\hat{\rho}_t)$ denotes the feedback controller as a function of the estimated state $\hat{\rho}_t$,
\item $\omega\geq0$ is the difference between the energies of the excited state and the ground state, $\eta\in(0,1]$ describes the efficiency of the detector, and $M>0$ is the strength of the interaction between the system and the probe. The estimated parameters $\hat{\omega}\geq0$, $\hat{\eta}\in(0,1]$ and $\hat{M}>0$, which may not equal to the actual ones.
\end{itemize}
By replacing $dY_t=dW_t+\sqrt{\eta M}\mathrm{Tr}(\sigma_z\rho_t)dt$ in the equation above, we obtain the following matrix-valued stochastic differential equations describing the time evolution of the pair $(\rho_t,\hat{\rho}_t)\in\mathcal{S}_2\times\mathcal{S}_2$,
 \begin{align}
 \setlength{\abovedisplayskip}{3pt}
\setlength{\belowdisplayskip}{3pt}
 d\rho_t&=L_{\omega, M}^u(\rho_t)dt+G_{\eta,M}(\rho_t)dW_t,\label{2D SME}\\
 d\hat{\rho}_t&=L_{\hat\omega, \hat M}^u(\hat\rho_t)dt +G_{ \hat \eta,\hat M}(\hat{\rho}_t)dW_t\nonumber\\
& +G_{\hat{\eta},\hat M}(\hat{\rho}_t)\big(\sqrt{\eta M}\mathrm{Tr}(\sigma_z\rho_t)-\sqrt{\hat{\eta} \hat{M}}\mathrm{Tr}(\sigma_z\hat{\rho}_t)\big)dt\label{2D SME filter}
 \end{align}
If $u\in\mathcal{C}^1(\mathcal{S}_2,\mathbb{R})$, the existence and uniqueness of the solution of~\eqref{2D SME}--\eqref{2D SME filter} can be 
proved along the same lines of~\cite[Proposition~3.5]{mirrahimi2007stabilizing}. Similarly, it can be shown as in~\cite[Proposition~3.7]{mirrahimi2007stabilizing} that $(\rho_t,\hat\rho_t)$ is a strong Markov process in $\mathcal{S}_2\times \mathcal{S}_2$.

Recall that a density operator $\rho\in \mathcal{S}_2$ can be uniquely characterized by the Bloch sphere coordinates $(x,y,z)$ as 
\begin{equation*}
\setlength{\abovedisplayskip}{3pt}
\setlength{\belowdisplayskip}{3pt}
\rho=\frac{\mathds{1}+x\sigma_x+y\sigma_y+z\sigma_z}{2}=\frac 12
\begin{bmatrix}
1+z & x-iy\\
x+iy & 1-z
\end{bmatrix},
\end{equation*}
where the vector $(x,y,z)$ belongs to the ball
\begin{equation*}
\setlength{\abovedisplayskip}{3pt}
\setlength{\belowdisplayskip}{3pt}
\mathcal{B}:=\{(x,y,z)\in\mathbb{R}^3|\, x^2+y^2+z^2 \leq 1 \}.
\end{equation*}
The stochastic differential equations~\eqref{2D SME}--\eqref{2D SME filter} expressed in Bloch sphere coordinates take the following form
\begin{subequations}
\begin{align}
dx_t&=\left(\!-\omega y_t-\frac{M}{2}x_t+uz_t \!\right)dt\!-\!\sqrt{\eta M}x_tz_tdW_t,\label{SME Bloch:x}\\
dy_t&=\left(\omega x_t-\frac{M}{2} y_t \right)dt-\sqrt{\eta M}y_tz_tdW_t, \label{SME Bloch:y}\\
dz_t&=-ux_tdt+\sqrt{\eta M} ( 1-z^2_t )dW_t.\label{SME Bloch:z}
\end{align}
\label{SME Bloch}
\end{subequations}
\vspace{-20pt}
\begin{subequations}
\begin{align}
d\hat{x}_t\!\!=&\big(\!\!-\!\hat\omega \hat{y}_t\!-\!\frac{\hat M}{2}\hat{x}_t\!+\!u\hat{z}_t\!+\!\hat{x}_t\hat{z}_t\mathcal E(z_t,\hat{z}_t)\big)dt\!-\!\!\sqrt{\hat{\eta} \hat{M}}\hat{x}_t\hat{z}_tdW_t\label{SME filter Bloch:x}\\
d\hat{y}_t\!=&\big(\hat\omega \hat{x}_t-\!\frac{\hat M}{2} \hat{y}_t +\hat{y}_t\hat{z}_t\mathcal E(z_t,\hat{z}_t)\big)dt-\!\sqrt{\hat{\eta} \hat{M}}\hat{y}_t\hat{z}_tdW_t, \label{SME filter Bloch:y}\\
d\hat{z}_t\!=&\big(\!-\!u\hat{x}_t\!-\!(1-\hat{z}^2_t)\mathcal E(z_t,\hat{z}_t)\big)dt+\!\sqrt{\hat{\eta} \hat{M}}( 1-\hat{z}^2_t)dW_t.\label{SME filter Bloch:z}
\end{align}
\label{SME filter Bloch}
\end{subequations}
where $\mathcal E(z,\hat{z}):= \sqrt{\hat{\eta} \hat{M}}\big(\sqrt{\hat{\eta} \hat{M}}\hat{z}-\sqrt{\eta M}z\big)).$
\section{Basic stochastic tools}
In this section, we introduce some basic definitions and classical results which are fundamental for the rest of the paper.

Given a stochastic differential equation $dq_t=f(q_t)dt+g(q_t)dW_t$, where $q_t$ takes values in $Q\subset \mathbb{R}^p,$ the infinitesimal generator  is the operator $\mathscr{L}$ acting on twice continuously differentiable functions $V: Q \times \mathbb{R}_+ \rightarrow \mathbb{R}$ in the following way
\begin{align*}
\setlength{\abovedisplayskip}{3pt}
\setlength{\belowdisplayskip}{3pt}
\mathscr{L}V(q,t):=&\frac{\partial V(q,t)}{\partial t}+\sum_{i=1}^p\frac{\partial V(q,t)}{\partial q_i}f_i(q)\\
&+\frac12 \sum_{i,j=1}^p\frac{\partial^2 V(q,t)}{\partial q_i\partial q_j}g_i(q)g_j(q).
\end{align*}
It\^o's formula describes the variation of the function $V$ along solutions of the stochastic differential equation and is given as follows
\begin{equation*}
\setlength{\abovedisplayskip}{3pt}
\setlength{\belowdisplayskip}{3pt}
dV(q,t) = \mathscr{L}V(q,t)dt+\sum_{i=1}^p\frac{\partial V(q,t)}{\partial q_i}g_i(q)dW_t.
\end{equation*}
From now on, the operator $\mathscr{L}$ is associated with~\eqref{2D SME}--\eqref{2D SME filter}.

We recall that the Bures metric for the 2-level case, which measures the ``distance" between two density matrices $\rho^{(1)},\rho^{(2)}$ in $\mathcal{S}_2,$ is given by $$
\setlength{\abovedisplayskip}{3pt}
\setlength{\belowdisplayskip}{3pt}
d_B(\rho^{(1)},\rho^{(2)}):=\sqrt{2-2\sqrt{\mathcal{F}(\rho^{(1)},\rho^{(2)})}},
$$
where 
$
\mathcal{F}(\rho^{(1)},\rho^{(2)})\!:=\!\mathrm{Tr}(\rho^{(1)},\rho^{(2)})\!+\!2\sqrt{\det(\rho^{(1)})\det(\rho^{(2)})}.
$
In particular, the Bures distance between $\rho\in\mathcal{S}_2$ and a pure state $\boldsymbol\rho=\psi\psi^*$ with $\psi\in\mathbb{C}^2$, is given by
$
d_B(\rho,\boldsymbol \rho)=\sqrt{2-2\sqrt{\psi^*\rho\psi}}.
$
In view of defining the notion of stochastic exponential  stability for the coupled system~\eqref{2D SME}--\eqref{2D SME filter}, we introduce the distance 
$$
\mathbf{d}_B\big((\rho^{(1)},\hat\rho^{(1)}),(\rho^{(2)},\hat \rho^{(2)})\big):=d_B(\rho^{(1)},\rho^{(2)})+d_B(\hat{\rho}^{(1)},\hat{\rho}^{(2)})
$$
between two elements of $\mathcal{S}_2\times\mathcal{S}_2$.
%
%
%
We denote the ball of radius $r$ around $(\rho,\hat \rho)$ as
$$
\mathbf{B}_r(\rho,\hat \rho) := \{(\sigma,\hat{\sigma}) \in\mathcal{S}_2\times\mathcal{S}_2|\, \mathbf{d}_B\big((\rho,\hat{\rho}),(\sigma,\hat \sigma)\big) < r\}.
$$
\begin{definition}
An equilibrium $(\boldsymbol \rho,\hat{\boldsymbol\rho})$ of the coupled system~\eqref{2D SME}--\eqref{2D SME filter} is said to be
\emph{almost surely exponentially stable} if
\begin{equation*}
\setlength{\abovedisplayskip}{3pt}
\setlength{\belowdisplayskip}{3pt}
\limsup_{t \rightarrow \infty} \frac{1}{t} \log \mathbf {d}_B((\rho_t,\hat\rho_t),(\boldsymbol \rho,\hat{\boldsymbol\rho})) < 0, \quad a.s.
\end{equation*}
whenever $( \rho_0,\hat\rho_0) \in \mathcal{S}_2\times \mathcal{S}_2$. The left-hand side of the above inequality is called the \emph{sample Lyapunov exponent} of the solution. 
\end{definition}
Denote $\boldsymbol \rho_g:=\mathrm{diag}(1,0)$ and $\boldsymbol \rho_e:=\mathrm{diag}(0,1)$, which are the pure states corresponding to the eigenvectors of~$\sigma_z$. 
 Note that a pair $(\boldsymbol \rho,\hat{\boldsymbol\rho})$ is an equilibrium of~\eqref{2D SME}--\eqref{2D SME filter} if and only if $\{\boldsymbol \rho,\hat{\boldsymbol\rho}\}\subset \{\boldsymbol \rho_e,\boldsymbol\rho_g\}$ and $u(\hat{\boldsymbol\rho})=0$. 
In order to introduce the final result of this section, we recall that any stochastic differential equation in It\^o form in $\mathbb R^K$
\begin{equation*}
\setlength{\abovedisplayskip}{3pt}
\setlength{\belowdisplayskip}{3pt}
dx_t=\widehat X_0(x_t)dt+\sum^n_{k=1}\widehat X_k(x_t)dW^k_t, \quad x_0 = x,
\end{equation*}
can be written in the following Stratonovich form~\cite{rogers2000diffusions2}
\begin{equation*}
\setlength{\abovedisplayskip}{3pt}
\setlength{\belowdisplayskip}{3pt}
dx_t = X_0(x_t)dt+\sum^n_{k=1}X_k(x_t) \circ dW^k_t, \quad x_0 = x,
\end{equation*}
where 
$X_0(x)=\widehat X_0(x)-\frac{1}{2}\sum^K_{l=1}\sum^n_{k=1}\frac{\partial \widehat X_k}{\partial x_l}(x)(\widehat X_k)_l(x)$, $(\widehat X_k)_l$ denoting the component $l$ of the vector $\widehat X_k,$ and $X_k(x)=\widehat X_k(x)$ for $k\neq 0$.

The following classical theorem relates the solutions of a stochastic differential equation with those of an associated deterministic one.
 \begin{theorem}[Support theorem~\cite{stroock1972support}]
Let $X_0(t,x)$ be a bounded measurable function, uniformly Lipschitz continuous in $x$ and $X_k(t,x)$  be continuously differentiable in $t$ and twice continuously differentiable in $x$, with bounded derivatives, for $k\neq 0.$ Consider the Stratonovich equation
\begin{equation*}
\setlength{\abovedisplayskip}{3pt}
\setlength{\belowdisplayskip}{3pt}
dx_t = X_0(t,x_t)dt+\sum^n_{k=1}X_k(t,x_t) \circ dW^k_t, \!\quad x_0 = x.
\end{equation*}
Let $\mathbb{P}_x$ be the probability law of the solution $x_t$ starting at $x$. Consider in addition the associated deterministic control system
\begin{equation}
\setlength{\abovedisplayskip}{3pt}
\setlength{\belowdisplayskip}{3pt}
\frac{d}{dt}x_{v}(t) = X_0(t,x_{v}(t))+\sum^n_{k=1}X_k(t,x_{v}(t))v^k(t), \!\!\quad x_v(0)\! =\! x.
\label{det}
\end{equation}
with $v^k \in \mathcal{V}$, where $\mathcal{V}$ is the set of all piecewise constant functions from $\mathbb{R}_+$ to $\mathbb{R}$. Now we define $\mathcal{W}_x$ as the set of all continuous paths from $\mathbb{R}_+$ to $\mathbb R^K$ starting at $x$, equipped with the topology of uniform convergence on compact sets, and $\mathcal{I}_x$ as the smallest closed subset of $\mathcal{W}_x$ such that $\mathbb{P}_x(x_{\cdot} \in \mathcal{I}_x)=1$. Then,
$
\mathcal{I}_x = \overline{ \{ x_{v}(\cdot)\in\mathcal{W}_x|\, v \in \mathcal{V}^n\} } \subset \mathcal{W}_x.
$
\label{Support thm}
\end{theorem}
\section{Feedback exponential stabilization of quantum spin-$\frac{1}{2}$ systems}
Our aim here is to provide sufficient conditions on the feedback controller $u(\hat{\rho})$ and a valid domain of the estimated parameters $\hat{\omega}$, $\hat{M}$ and $\hat{\eta}$ allowing us to exponentially stabilize the coupled system~\eqref{2D SME}--\eqref{2D SME filter} towards the target state $(\boldsymbol\rho_e,\boldsymbol\rho_e)$. By symmetry, the case in which the target state is $(\boldsymbol\rho_g,\boldsymbol\rho_g)$ can be treated in the same manner.

By employing arguments similar to those in~\cite[Lemma 3.1]{liang2019exponential_TwoUnknown}, we obtain the following invariance properties for the coupled system~\eqref{2D SME}--\eqref{2D SME filter}.
\begin{lemma}
Let $(\rho_t,\hat\rho_t)$ be the solution of~\eqref{2D SME}--\eqref{2D SME filter} starting from $(\rho_0,\hat\rho_0)$. If $\rho_0>0$, then $\mathbb{P}(\rho_t>0,\,\forall t\geq0)=1$. Similarly, if $\hat{\rho}_0>0$, then $\mathbb{P}(\hat{\rho}_t>0,\,\forall t\geq0)=1$. In other words, the sets $\mathrm{int}(\mathcal{S}_2)\times\mathcal{S}_2$ and $\mathcal{S}_2\times\mathrm{int}(\mathcal{S}_2)$ are almost surely invariant for~\eqref{2D SME}--\eqref{2D SME filter}.
\label{Lemma:Invaraint}
\end{lemma}

We make the following hypothesis on the feedback controller.
\begin{itemize}
\item[\textbf{H}:] 
$u\in\mathcal{C}(\mathcal{S}_2,\mathbb{R})\cap\mathcal{C}^1(\mathcal{S}_2\setminus\{\boldsymbol\rho_e\},\mathbb{R})$, $u(\boldsymbol\rho_e)=0$ and $u(\boldsymbol\rho_g)\neq0$.
\end{itemize}
If \textbf{H} is satisfied, then the coupled system~\eqref{2D SME}--\eqref{2D SME filter} admits exactly  two equilibria : $(\boldsymbol\rho_e,\boldsymbol\rho_e)$ and $(\boldsymbol\rho_g,\boldsymbol\rho_e)$. 

The following result, analogous to~\cite[Lemma 3.3]{liang2019exponential_TwoUnknown}, provides sufficient conditions guaranteeing that $\rho_t$ and $\hat{\rho}_t$ immediately become positive definite,  almost surely.  
\begin{lemma}
Assume that $\eta,\hat{\eta}\in(0,1)$ and \textbf{H} is satisfied. Then, for all initial condition $(\rho_0,\hat{\rho}_0)\in \partial\mathcal({S}_2\times\mathcal{S}_2)\setminus\{(\boldsymbol\rho_e,\boldsymbol\rho_e)\cup(\boldsymbol\rho_g,\boldsymbol\rho_e)\}$,  $(\rho_t,\hat{\rho}_t)\in\mathrm{int}(\mathcal{S}_2)\times\mathrm{int}(\mathcal{S}_2)$ for all $t>0$ almost surely.
\label{Lemma:PosDef}
\end{lemma}

Next, we show the instability of the equilibrium $(\boldsymbol\rho_g,\boldsymbol\rho_e)$.
\begin{lemma}
Suppose that \textbf{H} is satisfied and $|u(\hat\rho)|\leq C(1-\mathrm{tr}(\hat\rho\boldsymbol\rho_e))^{\alpha}$ for $C>0$ and $\alpha>\frac12$, then there exists $\lambda>0$ such that, for all initial condition $(\rho_0,\hat{\rho}_0)\in\mathbf{B}_\lambda(\boldsymbol\rho_g,\boldsymbol\rho_e)\setminus(\boldsymbol\rho_g,\boldsymbol\rho_e)$, the trajectories of the coupled system~\eqref{2D SME}--\eqref{2D SME filter} exit $\mathbf{B}_\lambda(\boldsymbol\rho_g,\boldsymbol\rho_e)$ in finite time almost surely.
\label{Lemma:Unstability}
\end{lemma}
\proof
We first show that for a small enough neighborhood of $(\boldsymbol\rho_g,\boldsymbol\rho_e)$ there exists a constant $\Gamma>2\hat{\eta}\hat M$ such that $\mathscr{L}(1-\hat{z})\geq \Gamma(1-\hat{z})$.

Indeed, by using  the fact that $|\hat{x}|\leq\sqrt{1-\hat{z}}$, we get
\begin{equation*}
\setlength{\abovedisplayskip}{3pt}
\setlength{\belowdisplayskip}{3pt}
\begin{split}
\mathscr{L}(1-\hat{z})&=u\hat{x}+(1-\hat{z}^2)\mathcal E(z,\hat{z})\\
&\geq \Big[-C(1-\hat{z})^{\alpha-\frac12}+(1+\hat{z})\mathcal{E}(z,\hat{z})\Big](1-\hat{z}).
\end{split}
\end{equation*}
The bracketed expression converges to $2\hat{\eta}\hat M+2\sqrt{\eta\hat{\eta}M\hat M}$ as $ (z,\hat z)$ converges to $(-1,1)$, so that for any $\Gamma\in (2\hat{\eta}\hat M,2\hat{\eta}\hat M+2\sqrt{\eta\hat{\eta}M\hat M})$ there exists a small enough neighborhood $\mathcal{U}$ of $(\boldsymbol\rho_g,\boldsymbol\rho_e)$ such that $\mathscr{L}(1-\hat{z})\geq \Gamma(1-\hat{z})$ holds true.
Let $\tau$ be the first exit time from $\mathcal{U}$.
Due to Lemma~\ref{Lemma:Invaraint}, we can apply It\^o's formula to $\log (1-\hat{z})$, obtaining $\mathscr{L}\log(1-\hat{z})\geq \Gamma-2\hat{\eta}\hat M>0$ on $\mathcal{U}$. By applying Dynkin formula~\cite{rogers2000diffusions2} to $\log (1-\hat{z})$ we obtain the following 
\begin{equation*}
\setlength{\abovedisplayskip}{3pt}
\setlength{\belowdisplayskip}{3pt}
\begin{split}
(\Gamma-2\hat{\eta}\hat M)\mathbb{E}(\tau)&\leq \mathbb{E}(\log(1-\hat{z}_{\tau}))-\log(1-\hat{z}_0)\\
&\leq \log 2-\log(1-\hat{z}_0).
\end{split}
\end{equation*}
Then by Markov inequality, we get
\begin{equation*}
\setlength{\abovedisplayskip}{3pt}
\setlength{\belowdisplayskip}{3pt}
\mathbb{P}(\tau=\infty)=\lim_{m\rightarrow\infty}\mathbb{P}(\tau>m)\leq \lim_{m\rightarrow\infty}\mathbb{E}(\tau)/m=0.
\end{equation*}
The proof is complete.
\hfill$\square$

Denote by $\tau_r$ the first time such that the trajectories of the coupled system~\eqref{2D SME}--\eqref{2D SME filter}
enter inside $\mathbf{B}_r(\boldsymbol\rho_e,\boldsymbol\rho_e),$ that is 
\[\setlength{\abovedisplayskip}{3pt}
\setlength{\belowdisplayskip}{3pt}
\tau_r:=\inf\{t>0|\,(\rho_t,\hat{\rho}_t)\in\mathbf{B}_r(\boldsymbol\rho_e,\boldsymbol\rho_e)\}.\]
We have the following lemma.
\begin{lemma}
Consider the coupled system~\eqref{2D SME}--\eqref{2D SME filter} and 
suppose that the feedback controller satisfies the assumptions of Lemma~\ref{Lemma:Unstability}. Then, for all $r>0$ and any given initial state $(\rho_0,\hat{\rho}_0)\in(\mathcal{S}_2\times\mathcal{S}_2)\setminus(\boldsymbol\rho_g,\boldsymbol\rho_e)$, $\mathbb{P}(\tau_r<\infty)=1$.
\label{Lemma:Reachability}
\end{lemma}
\proof
The lemma holds trivially true for $(\rho_0,\hat{\rho}_0)\in\mathbf{B}_r(\boldsymbol\rho_e,\boldsymbol\rho_e)$, as in that case $\tau_r=0$. Let us suppose that $(\rho_0,\hat{\rho}_0)\in(\mathcal{S}_2\times\mathcal{S}_2)\setminus\mathbf{B}_r(\boldsymbol\rho_e,\boldsymbol\rho_e)$. 

Consider the deterministic control system~\eqref{det} associated with~\eqref{2D SME}--\eqref{2D SME filter}. 
Following the proof of~\cite[Lemma 4.1]{liang2018exponential}, we can easily show that, for every initial condition $(\rho_0,\hat{\rho}_0)$ and $\epsilon>0$, there exist $T\in(0,\infty)$ and a piecewise constant controller $v(t)$ such that the corresponding trajectory reaches $\mathbf{B}_r(\boldsymbol\rho_e,\boldsymbol\rho_e)\cup\mathbf{B}_{\epsilon}(\boldsymbol\rho_g,\boldsymbol\rho_e)$ by time $T$.
Due to Theorem~\ref{Support thm}, there exists 
$\zeta\in(0,1)$ such that $\mathbb{P}_{(\rho_0,\hat{\rho}_0)}(\mu_{r,\epsilon}<T)>\zeta$~\footnote{Recall that $\mathbb{P}_{(\rho_0,\hat{\rho}_0)}$ corresponds to the joint probability law of $(\rho_t,\hat\rho_t)$ starting at $(\rho_0,\hat\rho_0);$ the associated expectation is denoted by $\mathbb{E}_{(\rho_0,\hat\rho_0)}$.} , where $\mu_{r,\epsilon}:=\inf\{t>0|\,(\rho_t,\hat{\rho}_t)\in\mathbf{B}_r(\boldsymbol\rho_e,\boldsymbol\rho_e)\cup\mathbf{B}_{\epsilon}(\boldsymbol\rho_g,\boldsymbol\rho_e)\}$. By the compactness of $\mathbf{S}_{r,\epsilon}:=(\mathcal{S}_2\times\mathcal{S}_2)\setminus(\mathbf{B}_r(\boldsymbol\rho_e,\boldsymbol\rho_e)\cup\mathbf{B}_{\epsilon}(\boldsymbol\rho_g,\boldsymbol\rho_e))$ and the Feller continuity of $(\rho_t,\hat{\rho}_t)$, 
we have $\zeta\geq \zeta_0>0$ for $(\rho_0,\hat{\rho}_0)\in \mathbf{S}_{r,\epsilon}$, so that
$
\sup_{(\rho_0,\hat{\rho}_0)\in\mathbf{S}_{r,\epsilon}}\mathbb{P}_{(\rho_0,\hat{\rho}_0)}(\mu_{r,\epsilon}\geq T)\leq 1-\zeta_0<1.
$
By Dynkin inequality~\cite{dynkin1965markov}, 
\begin{align*}
\setlength{\abovedisplayskip}{3pt}
\setlength{\belowdisplayskip}{3pt}
\sup_{(\rho_0,\hat{\rho}_0)\in\mathbf{S}_{r,\epsilon}}&\mathbb{E}_{(\rho_0,\hat{\rho}_0)}(\mu_{r,\epsilon})\\
&\leq \frac{T}{1-\sup_{(\rho_0,\hat{\rho}_0)\in\mathbf{S}_{r,\epsilon}}\mathbb{P}_{(\rho_0,\hat{\rho}_0)}(\mu_{r,\epsilon}\geq T)}\\
&\leq \frac{T}{\zeta_0}<\infty.
\end{align*}
Then by Markov inequality, for all $(\rho_0,\hat{\rho}_0)\in\mathbf{S}_{r,\epsilon}$, 
$\mathbb{P}_{(\rho_0,\hat{\rho}_0)}(\mu_{r,\epsilon}<\infty)=1.$

Choose $\lambda$ as in Lemma~\ref{Lemma:Unstability} and take $0<\epsilon<\delta<\lambda.$ Due to Theorem~\ref{Support thm}, the Feller continuity of the solutions $(\rho_t,\hat{\rho}_t)$ and the compactness of $\mathbf{S}_{r,\delta}$ there exists 
$\kappa\in(0,1)$ such that $\mathbb{P}_{(\rho_0,\hat{\rho}_0)}(\tilde\mu_{\epsilon}<\tau_r)\leq\kappa<1$ for all $(\rho_0,\hat{\rho}_0)\in\mathbf{S}_{r,\delta},$ where $\tilde{\mu}_{\epsilon}:=\inf\{t>0|\,(\rho_t,\hat{\rho}_t)\in\mathbf{B}_{\epsilon}(\boldsymbol\rho_g,\boldsymbol\rho_e)\}$.

Now we define two sequences of stopping times $\{\nu_{\delta}^k\}_{k\geq1}$ and $\{\tilde{\mu}_{\epsilon}^k\}_{k\geq0}$ such that $\tilde{\mu}_{\epsilon}^0=0$,  
\begin{equation*}
\begin{split}
\nu_{\delta}^{k+1}&=\inf\{t>\tilde{\mu}_{\epsilon}^{k}|\,(\rho_t,\hat{\rho}_t)\notin\mathbf{B}_{\delta}(\boldsymbol\rho_g,\boldsymbol\rho_e)\}\\
\tilde{\mu}_{\epsilon}^{k+1}&=\inf\{t>\nu_{\delta}^{k+1}|\,(\rho_t,\hat{\rho}_t)\in\mathbf{B}_{\epsilon}(\boldsymbol\rho_g,\boldsymbol\rho_e)\}.
\end{split}
\end{equation*}
In the following we calculate $\mathbb{P}_{(\rho_0,\hat{\rho}_0)}(\tilde{\mu}^m_{\epsilon}<\tau_r),$ 
\begin{align*}
&\mathbb{P}_{(\rho_0,\hat{\rho}_0)}(\tilde{\mu}^m_{\epsilon}<\tau_r)\\
&=\mathbb{P}_{(\rho_0,\hat{\rho}_0)}(\nu_{\delta}^{1}<\tau_r,\,\tilde{\mu}^1_{\epsilon}<\tau_r,\nu_{\delta}^{2}<\tau_r,\dots,\tilde{\mu}^m_{\epsilon}<\tau_r)\\
&=\mathbb{P}_{(\rho_{\nu_{\delta}^{1}},\hat{\rho}_{\nu_{\delta}^{1}})}(\tilde{\mu}^1_{\epsilon}<\tau_r)\dots\mathbb{P}_{(\rho_{\nu_{\delta}^{m}},\hat{\rho}_{\nu_{\delta}^{m}})}(\tilde{\mu}^m_{\epsilon}<\tau_r)\leq \kappa^m.
\end{align*}
For the above calculations, we used the strong Markov property of $(\rho_t,\hat\rho_t)$ and Lemma~\ref{Lemma:Unstability}, i.e., for all $(\rho_0,\hat{\rho}_0)\in\mathbf{B}_{\delta}(\boldsymbol\rho_g,\boldsymbol\rho_e)\setminus(\boldsymbol\rho_g,\boldsymbol\rho_e)$, $\mathbb{P}_{(\rho_{\tilde{\mu}_{\epsilon}^k},\hat{\rho}_{\tilde{\mu}_{\epsilon}^k})}(\nu_{\delta}^{k+1}<\infty)=1$.

Thus, for all $(\rho_0,\hat{\rho}_0)\in(\mathcal{S}_2\times\mathcal{S}_2)\setminus(\boldsymbol\rho_g,\boldsymbol\rho_e)$, we have $\mathbb{P}_{(\rho_0,\hat{\rho}_0)}(\tau_r=\infty)=\mathbb{P}_{(\rho_0,\hat{\rho}_0)}(\tilde{\mu}^m_{\epsilon}<\infty,\,\forall m>0)=0$. 
Then, the proof is complete.
\hfill$\square$

The following result provides general Lyapunov-type conditions ensuring exponential stabilization towards the target state $(\boldsymbol\rho_e,\boldsymbol\rho_e)$.
\begin{theorem}
Assume that $(\rho_0,\hat{\rho}_0)\in(\mathcal{S}_2\times\mathcal{S}_2)\setminus(\boldsymbol\rho_g,\boldsymbol\rho_e)$ and the feedback controller satisfies the assumptions of Lemma~\ref{Lemma:Unstability}. Additionally, suppose that there exists a function $V(\rho,\hat{\rho})$ such that $V(\boldsymbol\rho_e,\boldsymbol\rho_e)=0$, $V$ is positive outside the equilibrium $(\boldsymbol\rho_e,\boldsymbol\rho_e)$, continuous on $\mathcal{S}_2\times\mathcal{S}_2$ and twice continuously differentiable on the set $\mathrm{int}(\mathcal{S}_2)\times\mathrm{int}(\mathcal{S}_2)$. Moreover, suppose that there exist positive constants $C$, $C_1$ and $C_2$ such that 
\begin{itemize}
\item[(i)]\!$C_1 \mathbf{d}_B\big(\!(\rho,\hat{\rho}),(\boldsymbol\rho_e,\boldsymbol\rho_e)\!\big) \!\!\leq\!\! V\!(\rho,\hat{\rho})\!\! \leq\! \!C_2 \mathbf{d}_B\big(\!(\rho,\hat{\rho}),(\boldsymbol\rho_e,\boldsymbol\rho_e)\!\big)$
for every $(\rho,\hat\rho)\in\mathcal{S}_2\times\mathcal{S}_2$,
\item[(ii)] $\limsup_{(\rho,\hat{\rho})\rightarrow(\boldsymbol\rho_e,\boldsymbol\rho_e)}\frac{\mathscr{L}V(\rho,\hat{\rho})}{V(\rho,\hat{\rho})}\leq-C$.
\end{itemize}
Then, $(\boldsymbol\rho_e,\boldsymbol\rho_e)$ is a.s. exponentially stable for the coupled system~\eqref{2D SME}--\eqref{2D SME filter} with sample Lyapunov exponent less than or equal to $-{C}-\frac{K}{2}$, where $K:=\liminf_{(\rho,\hat{\rho})\rightarrow(\boldsymbol\rho_e,\boldsymbol\rho_e)}\varphi^2(\rho,\hat{\rho})$ with $\varphi(\rho,\hat{\rho}):=\frac{\partial V(\rho,\hat{\rho})}{\partial \rho}\frac{G_{\eta,M}(\rho)}{V(\rho,\hat{\rho})}+\frac{\partial V(\rho,\hat{\rho})}{\partial \hat{\rho}}\frac{G_{\hat{\eta},\hat{M}}(\rho)}{V(\rho,\hat{\rho})}$.
\label{Theorem:General Exp stabilization}
\end{theorem}
\textit{Sketch of the proof.} To prove Theorem~\ref{Theorem:General Exp stabilization} one may follow the same steps as in~\cite[Theorem 6.2]{liang2019exponential}, making use of the preliminary lemmas stated above. In particular the presence of a  function $V$ satisfying $(i)$ and such that $\mathscr{L}V\leq 0$ may be used to prove that $(\boldsymbol\rho_e,\boldsymbol\rho_e)$ is a locally stable equilibrium in probability, that is, given $\epsilon>0$ there exists $\delta>0$ and $\eta>0$ small enough such that 
$\mathbb{P}(\rho_t\in\mathbf{B}_{\epsilon}(\boldsymbol\rho_e,\boldsymbol\rho_e),\,\forall t\geq0)\geq \eta$, provided that $\rho_0\in \mathbf{B}_{\delta}(\boldsymbol\rho_e,\boldsymbol\rho_e)$. This, together with Lemma~\ref{Lemma:Reachability} and the strong Markov property of $(\rho_t,\hat{\rho}_t)$, implies the almost sure convergence to the target equilibrium. Finally,  in view of Lemma~\ref{Lemma:PosDef}, the $\mathcal{C}^2$ regularity of the function $V$ in $\mathrm{int}(\mathcal{S}_2)\times\mathrm{int}(\mathcal{S}_2)$ and the condition $(ii)$ imply 
\begin{equation*}
\setlength{\abovedisplayskip}{3pt}
\setlength{\belowdisplayskip}{3pt}
\limsup_{t \rightarrow \infty} \frac{1}{t} \log V(\rho_t,\hat\rho_t) \leq -C-\frac{K}2, \quad a.s.
\end{equation*}
(see~\cite[Theorem 6.2]{liang2019exponential} for more details).
The result then follows from condition~$(i)$.\qed

Next, under an additional assumption on the physical parameters $\eta,\hat\eta,M,\hat M$, we show the stabilizability of~\eqref{2D SME}--\eqref{2D SME filter} by explicitly exhibiting a Lyapunov function satisfying the assumptions of Theorem~\ref{Theorem:General Exp stabilization}.
\begin{theorem}
Consider the coupled system~\eqref{2D SME}--\eqref{2D SME filter} with $(\rho_0,\hat{\rho}_0)\in(\mathcal{S}_2\times\mathcal{S}_2)\setminus(\boldsymbol\rho_g,\boldsymbol\rho_e)$. 
If $\hat{\eta}\hat{M}<4\eta M$ 
and 
$u$ satisfies the assumptions of Lemma~\ref{Lemma:Unstability}, 
then $(\boldsymbol\rho_e,\boldsymbol\rho_e)$ is almost surely exponentially stable with sample Lyapunov exponent less than or equal to $-\sqrt{\hat{\eta}\eta\hat{M}M}-\frac12\min \{\eta M-\hat\eta\hat M,0\}$.
\label{Theorem:Feedback Exp stabilization}
\end{theorem}
\proof
We set
$V(\rho,\hat{\rho})=\sqrt{1-z}+\sqrt{1-\hat{z}}$
as a candidate Lyapunov function, and we show that it satisfies the assumptions of Theorem~\ref{Theorem:General Exp stabilization}. 

The nonnegative function $V$ is equal to zero at the equilibrium, it is continuous on $\mathcal{S}_2\times\mathcal{S}_2$ and twice continuously differentiable on $\mathrm{int}(\mathcal{S}_2)\times\mathrm{int}(\mathcal{S}_2)$. 

Condition \emph{(i)} of Theorem~\ref{Theorem:General Exp stabilization} follows from straightforward computations. 

We show that the condition \emph{(ii)} holds true. 
The infinitesimal generator of the candidate Lyapunov function is given by
$
\mathscr{L}V(\rho,\hat{\rho})=uU_1(\rho,\hat{\rho})+U_2(\rho,\hat{\rho}),
$
where
\begin{align*}
\setlength{\abovedisplayskip}{3pt}
\setlength{\belowdisplayskip}{3pt}
U_1(\rho,\hat{\rho})&=\frac12\big(x(1-z)^{-\frac12}+\hat{x}(1-\hat{z})^{-\frac12}\big),\\
U_2(\rho,\hat{\rho})&=-\frac18 \Big[\eta M(1+z)^2\sqrt{1-z}+\hat{\eta}\hat{M}(1+\hat{z})^2\sqrt{1-\hat{z}}\Big]\\
&~~+\frac12\sqrt{\hat{\eta}\hat{M}}(1+\hat{z})\big(\sqrt{\hat{\eta}\hat{M}}\hat{z}-\sqrt{\eta M}z\big)\sqrt{1-\hat{z}}.
\end{align*}
 Using the fact that $|x|\leq \sqrt{2(1-z)}$ and $|\hat{x}|\leq\sqrt{2(1-\hat{z})}$ we get that $|U_1(\rho,\hat{\rho})|\leq \sqrt{2}$.
Since $|u|\leq c(1-\mathrm{tr}(\hat\rho\boldsymbol\rho_e))^{\alpha} = c(\frac{1-\hat{z}}2)^{\alpha}$ for some $c>0$ and $\alpha>\frac12$, we then have
 \[\setlength{\abovedisplayskip}{3pt}
\setlength{\belowdisplayskip}{3pt}
\lim_{(\rho,\hat{\rho})\rightarrow(\boldsymbol\rho_e,\boldsymbol\rho_e)}\frac{uU_1(\rho,\hat{\rho})}{V(\rho,\hat{\rho})} = 0.\]
Hence
\begin{align}
\setlength{\abovedisplayskip}{3pt}
\setlength{\belowdisplayskip}{3pt}
&\limsup_{(\rho,\hat{\rho})\rightarrow(\boldsymbol\rho_e,\boldsymbol\rho_e)}\frac{\mathscr{L}V(\rho,\hat{\rho})}{V(\rho,\hat{\rho})}=\limsup_{(\rho,\hat{\rho})\rightarrow(\boldsymbol\rho_e,\boldsymbol\rho_e)}\frac{U_2(\rho,\hat{\rho})}{V(\rho,\hat{\rho})}\nonumber\\
&= \limsup_{(z,\hat{z})\rightarrow(1,1)}\frac{-\eta M \sqrt{1-z} - (2\sqrt{\hat{\eta}\eta\hat{M}M}-\hat\eta\hat{M})\sqrt{1-\hat{z}} }{2V(z,\hat{z})}\nonumber\\
&= \frac12 \max \{-\eta M,-2\sqrt{\hat{\eta}\eta\hat{M}M}+\hat\eta\hat M\}\nonumber\\
&=-\sqrt{\hat{\eta}\eta\hat{M}M}+\frac12\hat\eta\hat M,
\label{eq:nuav}
\end{align}
which is negative under the assumptions of the theorem.
This proves the condition (ii) of Theorem~\ref{Theorem:General Exp stabilization}.

Furthermore, in the notations of Theorem~\ref{Theorem:General Exp stabilization}, we have 
$$\setlength{\abovedisplayskip}{3pt}
\setlength{\belowdisplayskip}{3pt}
K=\liminf_{(\rho,\hat{\rho})\rightarrow(\boldsymbol\rho_e,\boldsymbol\rho_e)}\varphi^2=\min\{\eta M,\hat{\eta}\hat{M}\},$$ and the sample Lyapunov exponent is less than or equal to $-\sqrt{\hat{\eta}\eta\hat{M}M}-\frac12\min \{\eta M-\hat\eta\hat M,0\}$.\qed

\medskip

Next, we give an example of feedback controller satisfying the assumptions of the theorem above.
\begin{proposition}
\label{Prop:example}
Consider the coupled system~\eqref{2D SME}--\eqref{2D SME filter} with $(\rho_0,\hat{\rho}_0)\in(\mathcal{S}_2\times\mathcal{S}_2)\setminus(\boldsymbol\rho_g,\boldsymbol\rho_e)$ and suppose $\hat{\eta}\hat{M}<4\eta M$. 
Define the feedback controller 
\begin{equation}
\setlength{\abovedisplayskip}{3pt}
\setlength{\belowdisplayskip}{3pt}
u(\hat{\rho})=\alpha\big(1-\mathrm{Tr}(\hat{\rho}\boldsymbol\rho_e)\big)^\beta,
\label{Eq:Feedback}
\end{equation}
where $\alpha>0$ and $\beta\geq 1$. Then, $(\boldsymbol\rho_e,\boldsymbol\rho_e)$ is almost surely exponentially stable with sample Lyapunov exponent less than or equal to $-\sqrt{\hat{\eta}\eta\hat{M}M}-\frac12\min \{\eta M-\hat\eta\hat M,0\}$.
\end{proposition}
Note that in  \cite[Conjecture~4.4]{liang2019exponential_TwoUnknown} we proposed candidate feedback laws in order to exponentially stabilize spin-$J$ systems in the case of unknown initial states. Proposition~\ref{Prop:example} provides a positive answer to such a conjecture assuming,  in addition to unknown initial states, unawareness of the physical parameters.
\begin{remark}
By a symmetric reasoning,  the feedback controller 
\begin{equation}
\setlength{\abovedisplayskip}{3pt}
\setlength{\belowdisplayskip}{3pt}
u(\hat{\rho})=\alpha\big(1-\mathrm{Tr}(\hat{\rho}\boldsymbol\rho_g)\big)^\beta,
\label{Eq:Feedbackg}
\end{equation}
with $\alpha>0$ and $\beta\geq 1,$ almost surely exponentially stabilizes the coupled system~\eqref{2D SME}--\eqref{2D SME filter} with $(\rho_0,\hat{\rho}_0)\in(\mathcal{S}_2\times\mathcal{S}_2)\setminus(\boldsymbol\rho_e,\boldsymbol\rho_g),$ towards $(\boldsymbol\rho_g,\boldsymbol\rho_g)$  with sample Lyapunov exponent less than or equal to $-\sqrt{\hat{\eta}\eta\hat{M}M}-\frac12\min \{\eta M-\hat\eta\hat M,0\}$.
\end{remark}
\section{Simulation}
\begin{figure}[thpb]
\centering
\includegraphics[width=9.5cm]{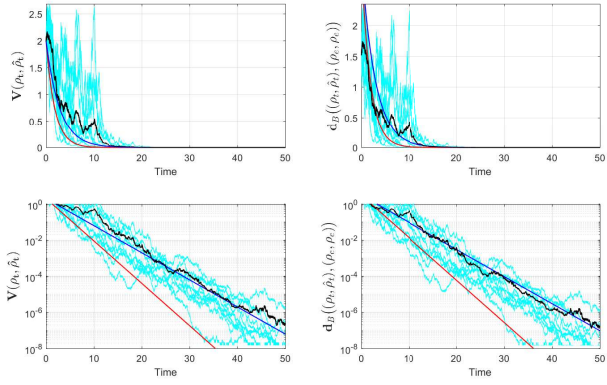}
\caption{Exponential stabilization of the coupled system towards $(\boldsymbol \rho_e,\boldsymbol \rho_e)$ with the feedback law~\eqref{Eq:Feedback} starting at $(x_0,y_0,z_0)=(1,0,0)$ and $(\hat{x}_0,\hat{y}_0,\hat{z}_0)=(0,1,0)$ with $\omega=0.3$, $\eta=0.3$, $M=1.3$, $\hat{\omega}=0.5$, $\hat{\eta}=0.5$, $\hat{M}=1.5$, $\alpha = 10$ and $\beta = 2$: the black curves represent the mean value of 10 arbitrary sample trajectories, the red curves represent the exponential reference with exponent $\nu_{\textrm{s}}=-0.5408$, the blue curve represents the exponential reference with exponent $\nu_{\textrm{av}}=-0.3458$. The figures at the bottom are the semi-log versions of the ones at the top.}
\label{RobustFilter_E}
\end{figure}

\begin{figure}[thpb]
\centering
\includegraphics[width=9.5cm]{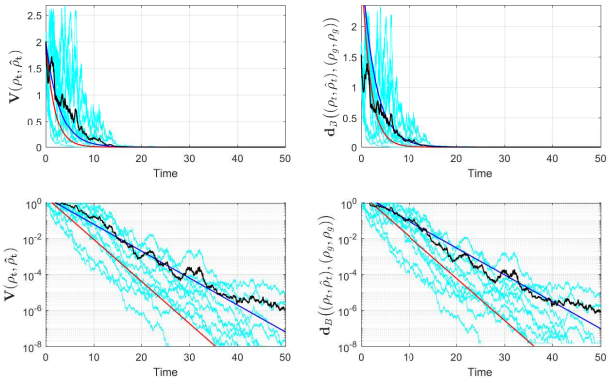}
\caption{Exponential stabilization of the coupled system towards $(\boldsymbol \rho_g,\boldsymbol \rho_g)$ with the feedback law~\eqref{Eq:Feedbackg} starting at $(x_0,y_0,z_0)=(1,0,0)$ and $(\hat{x}_0,\hat{y}_0,\hat{z}_0)=(0,1,0)$ with $\omega=0.3$, $\eta=0.3$, $M=1.3$, $\hat{\omega}=0.5$, $\hat{\eta}=0.5$, $\hat{M}=1.5$, $\alpha = 10$ and $\beta = 2$: the black curves represent the mean value of 10 arbitrary sample trajectories, the red curves represent the exponential reference with exponent $\nu_{\textrm{s}}=-0.5408$, the blue curves represent the exponential reference with exponent $\nu_{\textrm{av}}=-0.3458$. The figures at the bottom are the semi-log versions of the ones at the top.}
\label{RobustFilter_G}
\end{figure}

In this section, we first illustrate the convergence of the coupled system~\eqref{2D SME}--\eqref{2D SME filter}  starting at $(x_0,y_0,z_0)=(1,0,0)$ and $(\hat{x}_0,\hat{y}_0,\hat{z}_0)=(0,1,0)$ towards the target state $(\boldsymbol\rho_e,\boldsymbol\rho_e)$ by applying a feedback law of the form~\eqref{Eq:Feedback}. This is shown  in Figure~\ref{RobustFilter_E}. Then, in Figure~\ref{RobustFilter_G}, we show the convergence of the coupled system starting at the same initial states,  towards the target state $(\boldsymbol\rho_g,\boldsymbol\rho_g)$ by a feedback law of the form~\eqref{Eq:Feedbackg}. 

By Equation~\eqref{eq:nuav}, heuristically we have that the rate of convergence of the expectation of the Lyapunov function is less than or equal to $\nu_{\textrm{av}}:=-(\hat{\eta}\eta\hat{M}M)^{1/2}+\hat\eta\hat M/2$. This property is confirmed through simulations, see Fig.~\ref{RobustFilter_E} and Fig.~\ref{RobustFilter_G} (for the target state $(\boldsymbol\rho_g,\boldsymbol\rho_g),$ we take $V(\rho,\hat{\rho})=\sqrt{1+z}+\sqrt{1+\hat{z}}$). In the figures, the blue curves represent the exponential reference with the exponent $\nu_{\textrm{av}}$ and the black curves describe the mean values of the Lyapunov functions (Bures distances) of ten samples. On the figures, in particular in the semi-log versions, we can see that the black and the blue curves have similar asymptotic behaviors.  The red curves describe the exponential reference with exponent $\nu_{\textrm{s}}:=-(\hat{\eta}\eta\hat{M}M)^{1/2}-\frac12\min \{\eta M-\hat\eta\hat M,0\}$ and the cyan curves represent the behaviors of ten sample trajectories. We observe that the red cuves and the cyan curves have similar asymptotic behaviors.
\section{Conclusion}
In this paper, we studied the robustness of the stabilizing feedback strategy proposed in~\cite{liang2018exponential} for the case of spin-$\frac 12$ systems if initial states and physical parameters are unknown. We showed such a robustness property by analyzing the asymptotic behavior of the coupled system describing the evolutions of the quantum filter and the associated estimated state under appropriate assumptions on the feedback controller. More precisely, we showed exponential stabilization of the coupled system towards a pair $(\bar\rho,\bar\rho)$, with $\bar\rho$ being a chosen eigenstate of the measurement operator $\sigma_z$. Moreover, we gave an example of feedback  control law proving~\cite[Conjecture~4.4]{liang2019exponential_TwoUnknown} for spin-$\frac 1 2$ systems and supposing, in addition to unknown initial states and unlike~\cite{liang2019exponential_TwoUnknown}, that the exact values of the physical parameters are not accessible. A future research line 
will concern the robustness properties of the feedback controller considered in~\cite{liang2019exponential} for spin-$J$ systems.
\section{Acknowledgements}
This work is supported by the Agence Nationale de la Recherche projects Q-COAST ANR-19-CE48-0003 and QUACO ANR-17-CE40-0007.  The authors thank Pierre Rouchon for helpful discussions.

\bibliographystyle{plain}
\bibliography{Ref_Thesis_LIANG}

\end{document}